\documentclass[12pt]{article}
\usepackage[a4paper, left=.7in, right=.7in, top=1.5in, bottom=1.3in]{geometry}
\usepackage{amsmath, amssymb, graphicx, hyperref}

\newtheorem{theorem}{Theorem}[section]
\newtheorem{corollary}[theorem]{Corollary}
\newtheorem{lemma}[theorem]{Lemma}
\newtheorem{proposition}[theorem]{Proposition}

\newtheorem{definition}[theorem]{Definition}

\newtheorem{remark}[theorem]{Remark}
\newtheorem{example}[theorem]{Example}

\usepackage{fancyhdr}
\usepackage{makeidx}
\usepackage{hyperref}
\usepackage{moreverb}
\usepackage{tabularx}
\usepackage{tabulary}

\title{Softly $\pi g\widehat{D}$-Normal Spaces}
\author{
    Neeraj Kumar Tomar\\
    Gautam Buddha University, Greater Noida-201312, U.P. (India)\\
    Email: neer8393@gmail.com
    \and
    Fahed Zulfeqarr\\
    Gautam Buddha University, Greater Noida-201312, U.P. (India)\\
    Email: fahed@gbu.ac.in
    \and
    M. C. Sharma\\
    N. R. E. C. College, Khurja-203131, U.P. (India)\\
    Email: sharmamc2@gmail.com
}

\date{\today}
\date{}  

\begin{document}

\maketitle

\begin{abstract}
    The aim of this paper is to introduce a new class of softly normal called softly $\pi g\widehat{D}$ -normality by using $\pi g\widehat{D}$ -open sets and obtained several properties of such a space. We discuss many properties of this new space and we give some properties that connect this new spaces with some other topological spaces, also we present some examples and counter examples that show the relationships between softly $\pi g\widehat{D}$ -normal spaces and some other topological spaces, also we introduced the concept of $\pi g\widehat{D}$ -normal, almost $\pi g\widehat{D}$ -normal, quasi $\pi g\widehat{D}$ -normal, mildly $\pi g\widehat{D}$ -normal. The main result of this paper is that softly $\pi g\widehat{D}$ -normality is a topological property and it is a hereditary property with respect to closed domain subspaces. Moreover, we obtain some new characterizations and preservation theorems of softly $\pi g\widehat{D}$ -normal spaces. We insure that existence of utility for new results of softly $\pi g\widehat{D}$ -normality using separation axioms in topological spaces which is separate on a known separation axioms in topological spaces.
\end{abstract}


\section{Introduction}
\label{sec:intro}
In 1968, the notion of a quasi-normality is a weaker form of normality was introduced by Zaitsev\cite{ref1}. In 1970, the concept of almost normality was introduced by Singal and Arya\cite{ref2}. In 1973, the notion of mild normality was introduced by Shchepin\cite{ref3} and Singal and Singal\cite{ref4} independently. In 1990, Lal and Rahman\cite{ref5} have further studied notion of quasi normal and mildly normal spaces. In 2008, $\pi$-normal topological spaces were introduced by Kalantan\cite{ref6}. In 2011, Thabit and Kamarulhaili\cite{ref7} presented some characterizations of weakly (resp. almost) regular spaces. Also object of this paper is to present some conditions to assure that the product of two spaces will be $\pi$-normal. In 2012, Thabit and Kamarulhaili\cite{ref8} introduced a weaker version of $p$-normality called $\pi p$-normality and obtained some basic properties, examples, characterizations and preservation theorems of this property are presented. In 2014, Patil, Benchalli and Gonnagar\cite{ref9} introduced and studied two new classes of spaces, namely $\omega \alpha$-normal and $\omega\alpha$–regular spaces and obtained their properties by utilizing $\omega\alpha$-closed sets. In 2015, Hamant et al.\cite{ref10} introduced a new class of normal spaces called $\pi g\beta$-normal spaces, defined using $\pi g\beta$-open sets. They demonstrated that $\pi g\beta$-normality is a topological property and is hereditary with respect to $\pi$-open and $\pi g\beta$-closed subspaces. Furthermore, they established characterization and preservation theorems for $\pi g\beta$-normal spaces. In the same year, M. C. Sharma et al.\cite{ref11} introduced a weaker version of normality termed softly-normality. They proved that softly-normality is a property implied by both quasi-normality and almost normality, and investigated several properties of such spaces.\\
In this paper we present our work in five sections.\\
\textbf{Section 1.} First section have introduction and Literature reviews related to the research work.\\
\textbf{Section 2.} In this section, we define generalized $g\widehat{D}$-closed, $g\widehat{D}$-open, $\pi g\widehat{D}$-closed and $\pi g\widehat{D}$-open sets. We discuss some related examples and relationship with other class of closed sets.\\
\textbf{Section 3.} In this section, we defined softly $\pi g\widehat{D}$-normal, almost $\pi g\widehat{D}$-normal, quasi $\pi g\widehat{D}$-normal, $\pi g\widehat{D}$-normal and mildly $\pi g\widehat{D}$-normal spaces, we give related example and investigate the equivalent theorem related softly $\pi g\widehat{D}$-normal spaces.\\
\textbf{Section 4.} In this section, we investigate some function on softly $\pi g\widehat{D}$-normal spaces and investigate some theorems.\\
\textbf{Section 5.} In this section, we present some preservation theorem and other characterization of softly $\pi g\widehat{D}$-normal spaces.\\

\section{Preliminaries}
\label{sec:Preliminaries}
Throughout this paper, spaces $(X, \tau)$, $(Y, \sigma)$, and $(Z, \gamma)$ always mean topological spaces $X$, $Y$ and $Z$ on which no separation axioms are assumed unless explicitly stated. Let $A$ be a subset of a space $X$. The closure of $A$ and interior of $A$ are denoted by $cl(A)$ and $int(A)$ respectively.
\begin{definition}
    A subset $A$ of a topological space $X$ is said to be
    \begin{enumerate}
        \item \textbf{regular open}\cite{ref12} if $A$ = $int(cl(A))$.
        \item \textbf{pre-open}\cite{ref13} if $A \subset int(cl(A))$.
      \item \textbf{semi-open}\cite{ref14} if $A \subset cl(int(A))$.
    \item \textbf{$\alpha$-open} \cite{ref15} if $A \subset int(cl(int(A))).$
   \item \textbf{semi-preopen}\cite{ref16} if $A \subset cl(int(cl(A)))$.
   \item \textbf{ $\pi$-open} Finite union of regular open set is said to be $\pi$- open.
    \end{enumerate}   
    The complement of a regular-open (resp. pre-open, semi-open, $\alpha$-open, semi-preopen, $\pi$- open) set is called \textbf{regular-closed} (resp. \textbf{pre-closed, semi-closed, $\alpha$-closed, semi-preclosed, $\pi$- closed}) set.\\ 
    The intersection of all semi-preclosed (resp. pre-closed, semi-closed) sets containing $A$ is called the \textbf{semi-preclosure} (resp. \textbf{preclosure, semi-closure}) of $A$ and is denoted by \textbf{$spcl(A)$} (resp. \textbf{$pcl(A)$, $scl(A)$}). The \textbf{semi-preinterior} (resp. \textbf{preinterior, semi-interior}) of $A$, denoted by \textbf{$sp$-$int(A)$} (resp. \textbf{$p$-$int(A)$, $s$-$int(A)$}) is defined to be the union of all semi-preopen (resp. pre-open, semi-open) sets contained in $A$.
\end{definition}

\begin{definition}
    A subset $A$ of a topological spaces $X$ is called 
    \begin{enumerate}
        \item \textbf{generalized closed} (briefly \textbf{$g$-closed)}\cite{ref17} if $cl(A) \subset U$ whenever $A \subset U$ and $U$ is open in $X$.
        \item \textbf{$\pi g$-closed}\cite{ref18} if $cl(A) \subset U$  whenever  $A \subset U$  and  $U$  is $\pi$-open in  $X$.
        \item \textbf{$\alpha g$-closed}\cite{ref19} if  $\alpha cl(A) \subset U$  whenever $A \subset U$  and  $U$  is open in  $X$.
        \item \textbf{$\pi g\alpha$-closed}\cite{ref20} if  $\alpha cl(A) \subset U$  whenever  $A \subset U$  and  $U$  is $\pi$-open in  $X$.
        \item \textbf{generalized semi-preclosed} (briefly \textbf{$gsp$-closed)}\cite{ref21} if  $sp cl(A) \subset U$  whenever  $A \subset U$  and  $U$  is open in  $X$.
        \item \textbf{regular generalized closed} (briefly \textbf{$rg$-closed)}\cite{ref22} if  $ cl(A) \subset U$  whenever  $A \subset U$  and  $U$  is regular open in  $X$.
        \item \textbf{generalized pre-regular closed} (briefly \textbf{$gpr$-closed)}\cite{ref23} if  $p cl(A) \subset U$  whenever  $A \subset U$  and  $U$  is regular open in  $X$.
        \item \textbf{$w$-closed}\cite{ref24}, \textbf{$\widehat{g}$-closed} if $cl(A) \subset U$  whenever  $A \subset U$  and  $U$  is semi-open in  $X$.
        \item \textbf{pre-semi closed}\cite{ref25} if $spcl(A) \subset U$  whenever  $A \subset U$  and  $U$  is $g$-open in  $X$.
        \item \textbf{$D$-closed}\cite{ref26} if $pcl(A) \subset int(U)$  whenever  $A \subset U$  and  $U$  is $w$-open in  $X$.
    \end{enumerate}
    The complement of $g$- closed (resp. $\pi g$- closed, $\alpha g$- closed, $\pi g\alpha$- closed, $gsp$-closed, $rg$- closed, $gpr$- closed, $w$-closed, pre-semi-closed, $D$- closed) set is said to be \textbf{$g$- open} (resp. \textbf{$\pi g$- open, $\alpha g$- open, $\pi g\alpha$- open, $gsp$- open, $rg$- open, $gpr$- open, $w$- open, pre-semi- open, $D$- open}) set.
\end{definition}

\begin{definition}
    A subset $A$ of a topological space $X$ is called an $\widehat{D}$-closed\cite{ref27} set if $spcl(A) \subset U$ whenever $A \subset U$ and $U$ is $D$- open in $X$.\\
The class of all $\widehat{D}$-closed sets in $X$ is denoted by $\widehat{D}C(X)$.\\ 
That is $\widehat{D}C(X)$ = \{$A \subset X$ : $A$ is $\widehat{D}$-closed in $X$\}.\\
The complement of a $\widehat{D}$-closed set is called \textbf{$\widehat{D}$-open} set. The intersection of all $\widehat{D}$-closed set containing $A$ is called $\widehat{D}$-closure of $A$ and denoted by $\widehat{D}cl(A)$. The union of all $\widehat{D}$-open set containing A is called $\widehat{D}$-interior of $A$ and denoted by $\widehat{D}$-$int(A)$.

\end{definition}
\begin{definition}
    A subset $A$ of a topological space $X$ is called a $g\widehat{D}$-closed set if $\widehat{D}cl(A) \subset U$ whenever $A \subset U$ and $U$ is  open in $X$.
\end{definition}

\begin{definition}
    A subset $A$ of a topological space $X$ is called a $\pi g\widehat{D}$-closed set if $\widehat{D}cl(A) \subset U$ whenever $A \subset U$ and $U$ is  $\pi$-open in $X$.\vspace{.5cm}
    \\
    The complement of a $g\widehat{D}$-closed (resp. $\pi g\widehat{D}$-closed) set is called \textbf{$g\widehat{D}$-open} (resp. \textbf{$\pi g\widehat{D}$-open)} set, And the intersection of all $\pi g\widehat{D}$-closed set containing $A$ is called \textbf{$\pi g\widehat{D}$-closure} of $A$ and denoted by $\pi g\widehat{D}$-$cl(A)$. The union of all $\pi g\widehat{D}$-open set containing $A$ is called \textbf{$\pi g\widehat{D}$-interior} of $A$ and denoted by $\pi g\widehat{D}$-$int(A)$.
\end{definition}

\begin{remark}
    We have the following implications for the properties of subsets:

\[
\begin{array}{ccc}
\text{closed} & \longrightarrow & g\text{-closed} \longrightarrow \pi g\text{-closed} \\
\downarrow &  & \downarrow \quad \quad \quad \downarrow \\
\alpha\text{-closed} & \longrightarrow & \alpha g\text{-closed} \longrightarrow \pi g\alpha\text{-closed} \\
\downarrow &  & \downarrow \quad \quad \quad \downarrow \\
\widehat{D}\text{-closed} & \longrightarrow & g\widehat{D}\text{-closed} \longrightarrow \pi g\widehat{D}\text{-closed}
\end{array}
\]

Where none of the implications is reversible as can be seen from the following examples:
\end{remark}

\begin{example}
     Let $X$ be the real numbers with the cofinite topology. Let $A$ = $X - \{0\}$, then $A$ is a $\pi g$-closed set which is not $g$- closed.
 \end{example}   
\begin{example}
    Let $X$ = $\{a, b, c\}$ and $\tau$ = $\{X, \phi, \{a\}, \{b\}, \{a, b\}\}$. Then $\{a, b\}$ is $\pi g\widehat{D}$- closed but not $g\widehat{D}$-closed in $X$. 
\end{example}
\begin{proposition}
    Every closed (resp. $\alpha$-closed, pre-closed, semi-closed) set is $\pi g\widehat{D}$-closed set.
\end{proposition}
\textbf{Proof.}
    Let $A$ be any closed set. Let $A \subset U$ and $U$ is $\pi$-open set in $X$. Then    $cl(A) \subset U$. But $\widehat{D}cl(A) \subset cl(A) \subset U$. Thus $A$ is $\pi g\widehat{D}$-closed. The proof follows from the facts that \\
$\widehat{D}cl(A) \subset scl(A) \subset cl(A)$ and $\widehat{D}cl(A) \subset pcl(A) \subset \alpha cl(A) \subset cl(A)$.

\begin{remark}
    The converse of the above proposition need not be true as seen the following example. 
\end{remark}
\begin{example}
    Let $X$ = $\{a, b, c, d, e\}$ and $\tau$ = $\{\phi, \{a\}, \{a, b\}, \{c, d\}. \{a, c, d\}, \{a, b, c, d\}, X\}$. Then the set $A$ = $\{a, b, d\}$ is $\pi g\widehat{D}$- closed but not closed (resp. not pre-closed, not $\alpha$- closed, not semi-closed).
\end{example}
\begin{proposition}
    Every $w$-closed set is $\pi g\widehat{D}$-closed set.
\end{proposition}
\textbf{Proof.}
    Let A be $w$-closed set. Let $A \subset U$ and $U$ is $\pi$-open. Then $cl(A) \subset U$. Since every $w$-closed set is pre-closed and every pre-closed set is semi-pre-closed, $A$ is semi –pre- closed. Then $cl(A) \subset pcl(A) \subset wcl(A)$, since every closed is $w$-closed, $wcl(A) \subset cl(A)$.\\
Because $\widehat{D}cl(A) \subset pcl(A) \subset wcl(A) \subset cl(A) \subset U$. Hence A is $\pi g\widehat{D}$.

\begin{remark}
    The converse of the above proposition need not be true as seen the following example.
\end{remark}
\begin{example}
    Let $X$ = $\{a, b, c, d, e\}$ and $\tau$ = $\{\phi, \{a\}, \{b\}, \{a, b\}, X\}$. Then the set $A$ = $\{a\}$ is $\pi g\widehat{D}$- closed but not $w$-closed.
\end{example}
\begin{example}
    Let $X$ = $\{a, b, c\}$ and $\tau$ = $\{\phi, \{a\}, \{b\}, \{a, b\}, X\}$. Then 
    \begin{enumerate}
        \item closed sets in $X$ are $\phi$, $X$, $\{b, c\}$, $\{a, c\}$, $\{c\}$.
        \item $g$-closed sets in $X$ are $\phi$, $X$, $\{a\}$, $\{c\}$, $\{a, c\}$, $\{a, b, c\}$, $\{a, c, d\}$.
        \item $\widehat{D}$-closed sets in $X$ are $\phi$, $X$, $\{a\}$, $\{b\}$,$\{c\}$,$\{a, c\}$,$\{b, c\}$.
        \item $g \widehat{D}$-closed sets in $X$ are $\phi$, $X$, $\{a\}$, $\{b\}$,$\{c\}$,$\{a, c\}$,$\{b, c\}$.
        \item $\pi g\widehat{D}$-closed sets in $X$ are $\phi$, $X$, $\{a\}$, $\{b\}$,$\{c\}$,$\{a, c\}$,$\{b, c\}$.
    \end{enumerate}
\end{example}

\section{Softly $\pi g\widehat{D}$-Normal Spaces}
\begin{definition}
    A topological space $X$ is said to be \textbf{softly normal}\cite{ref11} (\textbf{softly $\pi g\widehat{D}$-normal}) if for any two disjoint subsets $A$ and $B$ of $X$, one of which is $\pi$-closed and other is regularly closed, there exist disjoint open ($\pi g\widehat{D}$-open) sets $U$ and $V$ of $X$ such that $A \subset U$ and $B \subset V$.
\end{definition}
\begin{definition}
    A topological space $X$ is said to be \textbf{almost normal} (\textbf{almost $\pi g \widehat{D}$-normal}) if for every pair of disjoint subsets $A$ and $B$ of $X$, one of which closed and other is regularly closed, there exist disjoint open ($\pi g\widehat{D}$- open) sets $U$ and $V$ of $X$ such that $A \subset U$ and $B \subset V$.
\end{definition}

\begin{definition}
    A topological space $X$ is said to be \textbf{$\pi$- normal}\cite{ref3} (\textbf{$\pi g\widehat{D}$- normal}) if for any two disjoint closed subsets $A$ and $B$ of $X$, one of which is $\pi$-closed, there exist disjoint open ($\pi g\widehat{D}$-open) sets $U$ and $V$ of $X$ such that $A \subset U$ and $B \subset V$. 

\end{definition}
\begin{definition}
    A topological space $X$ is said to be \textbf{quasi-normal}\cite{ref1} (\textbf{quasi $\pi g\widehat{D}$-normal}) if for any two disjoint $\pi$-closed subsets $A$ and $B$ of $X$, there exist disjoint open ($\pi g\widehat{D}$-open) sets $U$ and $V$ of $X$ such that $A \subset U$ and $B \subset V$.
\end{definition}
\begin{definition}
    A topological space $X$ is said to be \textbf{mildly normal}\cite{ref13} (\textbf{mildly $\pi g\widehat{D}$-normal}) if for any two disjoint regularly closed subsets $A$ and $B$ of $X$, there exist disjoint open ($\pi g\widehat{D}$-open) sets $U$ and $V$ of $X$ such that $A \subset U$ and $B \subset V$.
\end{definition}

\begin{remark}
    By the definitions stated above we have the following diagram holds for a topological space $X$.
    \begin{flushleft}
        \vspace{4.0mm}
    quasi-normal $\rightarrow$ quasi-$\pi g \widehat{D}$-normal $\rightarrow$ softly-$\pi g \widehat{D}$-normal $\rightarrow$ mildly-$\pi g \widehat{D}$-normal\\
   \hspace{.4cm} $\uparrow$\\ normal $\rightarrow$ $\pi$-normal $\rightarrow$ almost-normal $\rightarrow$ softly-normal
$\rightarrow$ mildly-normal \\ \vspace{.4cm} \hspace{2.4cm}$\downarrow\hspace{2.8cm}\downarrow\hspace{3cm}\downarrow\hspace{3.9cm}\downarrow$ \\ \vspace{.4cm} \hspace{1.cm} $\pi g \widehat{D}$-normal $\rightarrow$ almost-$\pi g \widehat{D}$-normal $\rightarrow$ softly-$\pi g \widehat{D}$-normal $\rightarrow$ mildly-$\pi g \widehat{D}$-normal
    \end{flushleft}
Where none of the implications is reversible as can be seen from the following examples.
\end{remark}
\begin{example}
Let $X$ = $\{a, b, c, d\}$ and $\tau$ = $\{\phi, \{a\}, \{c\}, \{a, c\}, \{b, d\}, \{a, b, d\}, \{b, c, d\}, X\}$. The pair of disjoint $\pi$- closed subset of $X$ are $A$ = $\{a\}$ and $B$ = $\{c\}$. Also $U$ = $\{a\}$ and $V$ = $\{b, c, d\}$ are disjoint open sets such that $A \subset U$ and $B \subset V$. Hence $X$ is quasi- normal as well as quasi $\pi g\widehat{D}$- normal as well as softly $\pi g\widehat{D}$- normal and because every open set is $\pi g\widehat{D}$- open set. 
\end{example}
\begin{example}
    Let $X$ = $\{a, b, c, d\}$ and $tau$ = $\{\phi, \{a\}, \{b\}, \{a, b\}, \{c, d\}, \{a, c, d\}, \{b, c, d\}, X\}$. Then $A$ = $\{b\}$ is closed set and $B$ = $\{a\}$ is regularly closed sets there exist disjoint open sets\\ $U$ = $\{b, c, d\}$ and $V$ = $\{a\}$ of $X$ such that $A \subset U$ and      
 $B \subset V$. Hence $X$ is almost- normal as well as almost $\pi g\widehat{D}$- normal as well as softly $\pi g\widehat{D}$- normal and because every open set is $\pi g\widehat{D}$- open set. 
\end{example}
\begin{example}
    Let $X$ = $\{a, b, c\}$ and $\tau$ = $\{\phi, \{a\}, \{a, b\}, \{a, c\}, X\}$. Then $X$ is almost- normal as well as almost $\pi g\widehat{D}$- normal, but it is not $\pi g\widehat{D}$- normal, since the pair of disjoint closed sets $\{b\}$ and $\{c\}$ have no disjoint $\pi g\widehat{D}$- open sets containing them. But it is not normal.
\end{example}
\begin{example}
Let $X$ = $\{a, b, c, d\}$ and $\tau$ = $\{\phi, \{a\}, \{b\}, \{a, b\}, \{c, d\}, \{a, c, d\}, \{b, c, d\}, X\}$. Then $A$ = $\{b\}$ is closed set and $B$ = $\{a\}$ is regularly closed sets there exist disjoint open sets $U$ = $\{b, c, d\}$ and $V$ = $\{a\}$ of $X$ such that $A \subset U$ and      
$B \subset V$. Hence $X$ is almost- normal as well as almost $\pi g\widehat{D}$- normal as well as softly $\pi g\widehat{D}$- normal and because every open set is $\pi g\widehat{D}$- open set. 
\end{example}
\begin{example}
Let $X$ = $\{a, b, c, d\}$ and $\tau$ = $\{\phi,\{a\},\{b\},\{d\},\{a,b\},\{b,d\},\{a,b,c\},\{a,b,d\},X\}$. Then $X$ is $\pi g
\widehat{D}$-normal.
\end{example}
\begin{theorem}
    For a topological space $X$ the following are equivalent:\\
    \textbf{a.} $X$ is softly $\pi g
\widehat{D}$-normal\\
    \textbf{b.} For every $\pi$-closed set $A$ and every regularly open set $B$ with $A\subset B$, there exists a $\pi g
\widehat{D}$-open set $U$ such that $A\subset U\subset \pi g\widehat{D}cl(U)\subset B$.   \\
    \textbf{c.} For every regularly closed set $A$ and every $\pi$- open set $B$ with $A \subset B$, there exists a $\pi g\widehat{D}$-open set $U$ such that
 $A\subset U\subset \pi g\widehat{D}cl(U)\subset B$.\\
 \textbf{d.} For every pair consisting of disjoint sets $A$ and $B$, one of which is $\pi$-closed and the other is regularly closed, there exist disjoint $\pi g\widehat{D}$-open sets $U$ and $V$ such that   $A  \subset U$, $B \subset V$ and $\pi g \widehat{D}$-$cl(U) \cap \pi g\widehat{D}$-$cl(V)$ = $\phi$.
\end{theorem}
\textbf{Proof.}
    $a \Rightarrow b$: Let $A$ be any $\pi$-closed set and $B$ be any regularly open set such that $A \subset B$. Then $A \cap (X - B)$ = $\phi$, where $(X - B)$ is regularly closed. Then there exist disjoint $\pi g \widehat{D}$-open sets $U$ and $V$ such that $A \subset U$ and 
    $(X - B) \subset V$. Since $U \cap V$ = $\phi$, then $\pi g \widehat{D}cl(U) \cap V$ = $\phi$. Thus
 $\pi g \widehat{D}cl(U) \subset (X - V) \subset (X - (X - B))$ = $B$. Therefore, 
$A \subset U \subset \pi g \widehat{D}cl(U) \subset B$.\\
$b \Rightarrow c$: Let $A$ be any regularly closed set and $B$ be any $\pi$-open set such that   $A \subset B$.\\ Then, $(X - B) \subset (X - A)$, where $(X - B)$ is $\pi$-closed and $(X - A)$ is regularly open. Thus by \textbf{b.,} there exists a $\pi g \widehat{D}$-open set $W$ such that
 $(X - B) \subset W \subset \pi g \widehat{D}cl (W) \subset (X - A)$.\\ 
Thus $A \subset (X - \pi g \widehat{D}cl (W)) \subset (X - W) \subset B$. So, we let $U$ = $(X - \pi g \widehat{D}cl(W))$, which is $\pi g \widehat{D}$-open and since  
$W \subset \pi g \widehat{D}cl(W)$, then $(X - \pi g \widehat{D}cl(W) \subset (X - W)$. Thus $U \subset (X - W)$, hence 
$\pi g \widehat{D}cl(U) \subset \pi g \widehat{D}cl(X - W)$ = $(X - W) \subset B$.\\
$c \Rightarrow d$: Let $A$ be any regular closed set and $B$ be any $\pi$-closed set with 
$A \cap B$ = $\phi$. Then $A \subset (X - B)$, where $(X - B)$ is $\pi$-open. By \textbf{c.,} there exists a $\pi g \widehat{D}$-open set $U$ such that,\\ $A \subset U \subset \pi g \widehat{D}cl (U) \subset (X - B)$. Now, $\pi g \widehat{D}cl (U)$ is $\pi g \widehat{D}$-closed. Applying \textbf{c.,} again we get a $\pi g \widehat{D}$-open set $W$ such that 
$A \subset U \subset \pi g \widehat{D}cl(U) \subset W \subset \pi g \widehat{D}cl(W) \subset (X - B)$.
 Let $V$ = $(X - \pi g \widehat{D}cl (W))$, then $V$ is $\pi g \widehat{D}$-open set and $B \subset V$. We have
 $(X - \pi g \widehat{D}cl (W)) \subset (X - W)$, hence $V \subset (X - W)$, thus
 $\pi g \widehat{D}cl (V) \subset \pi g \widehat{D}cl(X - W)$ = $(X - W)$. So, we have $\pi g \widehat{D}cl(U) \subset W$ and $\pi g \widehat{D}cl(V) \subset (X-W)$. Therefore $\pi g \widehat{D}cl(U) \cap \pi g \widehat{D}cl(V)$ = $\phi$.\\
 $d \Rightarrow a$: is clear.

\begin{theorem}
    For a topological space $X$ the following are equivalent: \\
    \textbf{a.} $X$ is softly $\pi g\widehat{D}$-normal,\\
    \textbf{b.} For every pair of sets $U$ and $V$, one of which $\pi$-open and other is regular open whose union is $X$, there exists disjoint $\pi g\widehat{D}$- closed sets $G$ and $H$ such that G$ \subset U$, $H \subset V$ and $G \cup H$ = $X$,\\
    \textbf{c.} For every $\pi$-closed set $A$ and every regular open set $B$ containing $A$, there exists a $\pi g\widehat{D}$-open set $V$ such that
    $A \subset   V \subset  \pi g\widehat{D}cl (V) \subset  B$.

\end{theorem}
\textbf{Proof.}
    $a \Rightarrow b$: Let $U$ be a $\pi$-open set and $V$ be a regular open set in a softly $\pi g\widehat{D}$-normal space $X$ such that $U \cup V$ = $X$. Then $(X - U)$ is $\pi$-closed set and $(X - V)$ is regular closed set with $(X - U) \cap (X - V)$ = $\phi$. By soft $\pi g\widehat{D}$-normality of $X$, there exist disjoint $\pi g\widehat{D}$-open sets $U_1$ and $V_1$ such that $X - U \subset U_1$ and $X - V \subset V_1$. Let $G$ = $X - U_1$ and $H$ = $X - V_1$. Then $G$ and $H$ are $\pi g\widehat{D}$-closed sets such that $G \subset U$, $H \subset V$ and $G \cup H$ = $X$.\\
    $b \Rightarrow c$: and $c \Rightarrow a$: are obvious.\\
    Using \textbf{Theorem 3.12.,} it is easy to show the following theorem, which is a Urysohn’s Lemma version for soft $\pi g\widehat{D}$-normality. A proof can be established by a similar way of the normal case.

\begin{theorem}
    A space $X$ is softly $\pi g\widehat{D}$-normal  if and only if for every pair of disjoint closed sets $A$ and $B$, one of which is $\pi$-closed and other is regularly closed, there exists a continuous function $f$ on $X$ into $[0, 1]$, with its usual topology, such that $f(A)$ = $\{0\}$ and $f(B)$ = $\{1\}$. It is easy to see that the inverse image of a regularly closed set under an open continuous function is regularly closed and the inverse image of a $\pi$-closed set under an open continuous function $\pi$-closed. We will use that in the next theorem.
\end{theorem}
\begin{theorem}
    Let $X$ is a softly $\pi g\widehat{D}$-normal space and $f : X \rightarrow Y$ is an open continuous injective function. Then $f(X)$ is a softly $\pi g\widehat{D}$-normal space.
\end{theorem}
\textbf{Proof.}
    Let $A$ be any $\pi$-closed subset in $f(X)$ and let $B$ be any regularly closed subset in $f(X)$ such that $A \cap B$ = $\phi$. Then $f^{-1}(A)$ is a $\pi$-closed set in $X$, which is disjoint from the regularly closed set $f^{-1}(B)$. Since $X$ is softly $\pi g\widehat{D}$-normal, there are two disjoint open sets $U$ and $V$ such that $f^{-1}(A) \subset U$ and $f^{-1}(B) \subset V$. Since $f$ is one-one and open, result follows. 

\begin{corollary}
    Soft $\pi g\widehat{D}$-normality is a topological property. 
\end{corollary}
\begin{lemma}
    Let $M$ be a closed domain subspace of a space $X$. If $A$ is a $\pi g\widehat{D}$-open set in $X$, then $A \cap M$ is $\pi g\widehat{D}$-open set in $M$.
\end{lemma}
\begin{theorem}
    A closed domain subspace of a softly $\pi g\widehat{D}$-normal is softly $\pi g\widehat{D}$-normal.
\end{theorem}
\textbf{Proof.}
    Let $M$ be a closed domain subspace of a softly $\pi g\widehat{D}$-normal space $X$. Let $A$ and $B$ be any disjoint closed sets in $M$ such that $A$ is regularly closed and $B$ is $\pi$-closed. Then, $A$ and $B$ are disjoint closed sets in $X$ such that $A$ is regularly closed and $B$ is $\pi$-closed in $X$. By soft $\pi g\widehat{D}$-normality of $X$, there exist disjoint $\pi g\widehat{D}$-open sets $U$ and $V$ of $X$ such that $A \subset U$ and $B \subset V$. By the above \textbf{Lemma 3.17.,} we have $U \cap M$ and $V \cap M$ are disjoint $\pi g\widehat{D}$-open sets in $M$ such that 
    $A \subset U \cap M$ and $B \subset V \cap M$. Hence, $M$ is softly $\pi g\widehat{D}$-normal subspace. Since every closed and open (clopen) subset is a closed domain, then we have the following corollary. 

\begin{corollary}
    Soft $\pi g\widehat{D}$-normality is a hereditary with respect to clopen subspaces.
\end{corollary}
\begin{definition} A space $X$ is said to be
    \begin{enumerate}
        \item \textbf{$\pi g\widehat{D}$-$T_2$} if for any distinct pair of points $x$ and $y$ in $X$, there exists $\pi g\widehat{D}$-open sets $U$ and $V$ in $X$ containing $x$ and $y$, respectively such that $U \cap V$ = $\phi$.
        \item \textbf{$\pi g\widehat{D}$-$T_3$} if it is softly $\pi g\widehat{D}$-normal as well as $\pi g\widehat{D}$-$T_1$ space.
    \end{enumerate}
\end{definition}
\begin{theorem}
    Every $\pi g\widehat{D}$-$T_3$ space is a $\pi g\widehat{D}$-$T_2$ space.
\end{theorem}
\textbf{Proof.}
    Let $X$ be $\pi g\widehat{D}$-$T_3$, so it is both $\pi g\widehat{D}$-$T_1$ and softly $\pi g\widehat{D}$-normal. Also $X$ is $\pi g\widehat{D}$-$T_1 \Rightarrow$ every singleton subset $\{x\}$ of $X$ is an $\pi g\widehat{D}$-closed. Let $\{x\}$ be an $\pi g\widehat{D}$-closed subset of $X$ and $y\in X - \{x\}$. Then we have $x\neq y$ since $X$ is softly $\pi g\widehat{D}$-normal, there exist disjoint $\pi g\widehat{D}$-open sets $U$ and $V$ such that $\{x\} \subset U$, $y\subset V$, and such that $U\cap V$ = $\phi $ (or) $U$ and $V$ are disjoint $\pi g\widehat{D}$-open sets containing $x$ and $y$ respectively. Since $x$ and $y$ are arbitrary, for every pair of distinct points, there exist disjoint $\pi g\widehat{D}$-open sets. Hence $X$ is $\pi g\widehat{D}$-$T_2$ space.

\section{Some function on softly $\pi g\widehat{D}$-normal spaces}
\begin{definition}
    A function $F : X\longrightarrow Y$ is said to be
    \begin{enumerate}
        \item \textbf{$\pi g\widehat{D}$-closed} if for every closed set $F$ of $X$, $f(F)$ is $\pi g\widehat{D}$-closed in $Y$.
        \item \textbf{almost $\pi g\widehat{D}$-closed} if for every regular closed set $F$ of $X$, $f(F)$ is $\pi g\widehat{D}$-closed in $Y$.
        \item \textbf{$\pi$-continuous}\cite{ref19} if $f^{-1}(F)$ is $\pi$- closed in $X$ for every closed set $F$ of $Y$.
        \item \textbf{$\pi g\alpha$-continuous} if $f^{-1}(F)$ is $\pi g\alpha$-closed in $X$ for every closed set $F$ of $Y.$
        \item \textbf{$\pi g\widehat{D}$-continuous} if $f^{-1}(F)$ is $\pi g\widehat{D}$-closed in $X$ for every closed set $F$ of $Y$.
        \item \textbf{almost continuous}\cite{ref4} if $f^{-1}(F)$ is closed in $X$ for every regular closed set $F$ of $Y$.
        \item \textbf{almost $\pi$-continuous}\cite{ref18} if $f^{-1}(F)$ is $\pi$-closed in $X$ for every regular closed set $F$ of $Y$.
        \item \textbf{almost $\pi g\alpha$-continuous} if $f^{-1}(F)$ is $\pi g\alpha$-closed in $X$ for every regular closed set $F$ of $Y$.
        \item \textbf{almost $\pi g\widehat{D}$-continuous} if $f^{-1}(F)$ is $\pi g\widehat{D}$-closed in $X$ for every regular closed set $F$ of $Y$.
        \item \textbf{rc-preserving}\cite{ref28} if $f(F)$ is regular closed in $Y$ for every $F\in RC(X)$.
        \item \textbf{softly $\pi g\widehat{D}$-irresolute} if for each $x$ in $X$ and each $\pi g\widehat{D}$- neighbourhood of $V$ of $f(x)$, $\pi g\widehat{D}$-$cl(f^{-1}(V))$ is a $\pi g\widehat{D}$-neighbourhood of $x$.
    \end{enumerate}
    From the definitions stated above, we obtain the following diagram:

\vspace{0.3cm}

\[
\begin{array}{cccccc}
\text{closed} & \Rightarrow & \text{almost-closed} \\
\Downarrow && \Downarrow \\
\alpha\text{-closed} & \Rightarrow & \text{almost } \widehat{D}\text{-closed} \\
\Downarrow && \Downarrow \\
g\alpha\text{-closed} & \Rightarrow & \text{almost } g\widehat{D}\text{-closed} \\
\Downarrow && \Downarrow \\
\pi g\alpha\text{-closed} & \Rightarrow & \text{almost } \pi g \widehat{D}\text{-closed}
\end{array}
\]
\vspace{.3cm}
Moreover, by the following examples, we realize that none of the implications is reversible.    
\end{definition}

\begin{example}
    $X$ = $\{a,b,c,d\}$, $\tau$ = $\{\phi,X,\{c\},\{a, b, d\}\}$ and $\sigma$ = $\{\phi,\{a\},\{c, d\},\{a, c, d\},\\ \{d\}, \{a, d\}, X\}$.
    Let $f : (X, \tau) \longrightarrow (X, \sigma)$ be the identity function, then $f$ is $\pi g\alpha$-closed as well as $\pi g\widehat{D}$-closed but not $\pi g$-closed. Since $A$ = $\{c\}$ is not $\pi g$-closed in $(X, \sigma)$.
\end{example}
\begin{example}
    Let $X$ = $\{a, b, c, d\}$, $\tau$ = $\{\phi,X,\{c\},\{b, d\}, \{a,b,d\}, \{b, c, d\}, X \}$ and $\sigma$ = $\{\phi, X, \{a\},\\ \{c, d\}, \{a, c, d\}, \{d\}, \{a, d\} X\}$. Let $f : (X,\tau) \longrightarrow (X,\sigma)$ be the identity function. Then $f$ is almost $\pi g\alpha$-closed as well as almost $\pi g\widehat{D}$- closed but not $\pi g\widehat{D}$- closed. Since $A$ = $\{a\}$ is not $\pi g\widehat{D}$- closed.
\end{example}

\begin{theorem}
     If $f : X \longrightarrow Y$ is an almost $\pi$-continuous and $\pi g\widehat{D}$-closed function, then $f(A)$ is $\pi g\widehat{D}$-closed in $Y$ for every $\pi g\widehat{D}$-closed set $A$ of $X$.
\end{theorem}
\textbf{Proof.}
    Let $A$ be any $\pi g\widehat{D}$-closed set $A$ of $X$ and $V$ be any $\pi$- open set of $Y$ containing $f(A)$. Since $f$ is almost $\pi$-continuous, $f^{-1}(V)$ is $\pi$-open in $X$ and $A \subset f^{-1}(V)$. Therefore $\widehat{D}cl(A) \subset f^{-1}(V)$ and hence $f(\widehat{D}cl(A)) \subset V$. Since $f$ is $\pi g\widehat{D}$-closed, $f(\widehat{D}cl(A)$ is $\pi g\widehat{D}$-closed in $Y$. And hence we obtain $\widehat{D}cl(f(A)) \subset \widehat{D}cl(f(\widehat{D}cl(A)))\subset V$. Hence $f(A)$ is $\pi g\widehat{D}$-closed in $Y$.

\begin{theorem}
    A surjection $f : X \longrightarrow Y$ is almost $\pi g\widehat{D}$-closed if and only if for each subset $S$ of $Y$ and each $U \in RO(X)$ containing $f^{-1}(S)$ there exists a $\pi g\widehat{D}$-open set $V$ of $Y$ such that $S \subset V$ and $f^{-1}(V) \subset U$. 
\end{theorem}
\textbf{Proof.}
    \textbf{Necessity.} Suppose that $f$ is almost $\pi g\widehat{D}$-closed. Let $S$ be a subset of $Y$ and $U \in RO(X)$ containing $f^{-1}(S)$. 
If $V$ = $Y - f(X - U)$, then $V$ is a $\pi g\widehat{D}$-open set of $Y$ such that $S \subset V$ and $f^{-1}(V) \subset U$.\\
\textbf{Sufficiency.} Let $F$ be any regular closed set of $X$. Then $f^{-1}(Y - f(F)) \subset X - F$ and $X - F \in RO(X)$. There exists $\pi g\widehat{D}$-open set $V$ of $Y$ such that $Y - f(F) \subset V$ and $f^{-1}(V) \subset X - F$. Therefore, we have 
$f(F)\subset Y - V$ and $F \subset X - f^{-1}(V) \subset f^{-1} (Y - V)$. Hence we obtain $f(F)$ = $Y - V$ and $f(F)$ is $g\widehat{D}$-closed in $Y$ which shows that $f$ is almost $\pi g\widehat{D}$-closed.

\section{Preservation theorems of softly $\pi g\widehat{D}$-normal spaces}

\begin{theorem}
    If $f : X \longrightarrow Y$ is a continuous softly $\pi g\widehat{D}$-closed surjection and $X$ is softly $\pi g\widehat{D}$-normal, then $Y$ is softly $\pi g\widehat{D}$-normal. 
\end{theorem}
\textbf{Proof.}
    Let $M_1$ and $M_2$ be any disjoint closed sets of $Y$. Since $f$ is continuous, $f^{-1}(M_1)$ and $f^{-1}(M_2)$ are disjoint closed sets of $X$. Since $X$ is softly $\pi g\widehat{D}$-normal, there exist disjoint\\ $U_1, U_2 \in \pi g\widehat{D}O(X)$ such that $f^{-1}(V_1)\subset U_i$ for $i = 1, 2$. Put $V_i$ = $Y - f(X - U_i)$, then $V_i$ is open in $Y$, $M_i\subset V_i$ and $f^{-1}(V_i)\subset U_i$ for $i = 1, 2$. Since $U_1\cap U_2$ = $\phi$ and $f$ is surjective;\\ we have $V_1\cap V_2$ = $\phi$. This shows that $Y$ is softly $\pi g\widehat{D}$-normal.

\begin{lemma}
    A subset $A$ of a topological space $X$ is $\pi g\widehat{D}$-open if and only if $F\subset \pi g\widehat{D}$-$int(A)$ whenever $F$ is closed and $F \subset A$.
\end{lemma}
\begin{corollary}
    If $f: X\longrightarrow Y$ is a closed $\pi g\widehat{D}$-irresolute injection and $Y$ is softly $\pi g\widehat{D}$-normal, then $X$ is softly $\pi g\widehat{D}$-normal. 
\end{corollary}
\begin{lemma}
    A function $f: X\longrightarrow Y$ is almost $\pi g\widehat{D}$-closed if and only if for each subset $B$ of $Y$ and each $U\in RO(X)$ containing $f^{-1}(B)$, there exists a $\pi g\widehat{D}$-open set $V$ of $Y$ such that $B \subset V$ and $f^{-1}(V)\subset U$. 
\end{lemma}
\begin{theorem}
    If $f : X \longrightarrow Y$ is an softly $\pi g\widehat{D}$-continuous, $rc$-preserving injection and Y is softly $\pi g\widehat{D}$-normal then $X$ is softly $\pi g\widehat{D}$-normal. 
\end{theorem}
\textbf{Proof.}
    Let $A$ and $B$ be any disjoint $\pi$-closed sets of $X$. Since $f$ is a $rc$-preserving injection, $f(A)$ and $f(B)$ are disjoint $\pi$-closed sets of $Y$. Since $Y$ is softly $\pi g\widehat{D}$-normal, there exist disjoint $\widehat{D}$-open sets $U$ and $V$ of $Y$ such that $f(A) \subset U$ and     $f(B) \subset V$. Now if $G$ = $int(cl(U))$ and $H$ = $int(cl(V))$. Then $G$ and $H$ are regular open sets such that $f(A) \subset G$ and $f(B) \subset H$. Since $f$ is almost $\pi g\widehat{D}$-continuous, $f^{-1}(G)$ and $f^{-1}(H)$ are disjoint $\pi g\widehat{D}$-open sets containing $A$ and $B$ which shows that $X$ is softly $\pi g\widehat{D}$-normal.

\begin{theorem}
    If $f: X\longrightarrow Y$ is $\pi$-continuous, $\pi g\widehat{D}$-closed surjection and X is softly $\pi g\widehat{D}$-normal space then $Y$ is softly $\pi g\widehat{D}$-normal.
\end{theorem}
\textbf{Proof.}
    Let $A$ and $B$ be any two disjoint closed sets of $Y$. Then $f^{-1}(A)$ and $f^{-1}(B)$ are disjoint $\pi$-closed sets of $X$. Since $X$ is softly $\pi g\widehat{D}$-normal, there exist disjoint softly $\pi g\widehat{D}$-open sets of $U$ and $V$ such that $f^{-1}(A)\subset U$ and $f^{-1}(B)\subset V$. Let $G$ = $int(cl(V))$ and $H$ = $int(cl(V))$. Then $G$ and $H$ are disjoint regular open sets of $X$ such that $f^{-1}(A)\subset G$ and $f^{-1}(B)\subset H$. Set $K$ = $Y - f(X - G)$ and $L$ = $Y - f (X - H)$. Then $K$ and $L$ are $\pi g\widehat{D}$-open sets of $Y$ such that $A \subset K$, $B \subset L$, $f^{-1}(K) \subset G$ , 
$f^{-1}(L) \subset H$. Since $G$ and $H$ are disjoint, $K$ and $L$ are disjoint. Since $K$ and $L$ are $\pi g\widehat{D}$-open and we obtain $A \subset\pi g\widehat{D}$-$int(K)$,
 $B\subset \pi g\widehat{D}$-$int(L)$ and $\pi g\widehat{D}$-$int(K)\cap \pi g\widehat{D}$-$int(L)$ = $\pi$. Therefore $Y$ is softly $\pi g\widehat{D}$-normal

\begin{theorem}
    Let $f : X\longrightarrow Y$ be an almost $\pi$-continuous and almost $\pi g\widehat{D}$-closed surjection. If $X$ is softly $\pi g\widehat{D}$-normal space then $Y$ is softly $\pi g\widehat{D}$-normal.
\end{theorem}
\textbf{Proof.}
    Let $A$ and $B$ be any disjoint $\pi$-closed sets of $Y$. Since $f$ is almost $\pi$-continuous, $f^{-1}(A)$, $f^{-1}(B)$ are disjoint closed subsets of $X$. Since $X$ is softly $\pi g\widehat{D}$-normal, there exist disjoint $\pi g\widehat{D}$-open sets $U$ and $V$ of $X$ such that $f^{-1}(A) \subset U$ and
$f^{-1}(B) \subset V$. Put $G$ = $int(cl(U))$ and $H$ = $int(cl(V))$. Then $G$ and $H$ are disjoint regular open sets of $X$ such that $f^{-1}(A)\subset G$ and $f^{-1}(B)\subset H$. By \textbf{Theorem 4.5.,} there exist $V$-open sets $K$ and $L$ of $Y$ such that $A \subset K$, $B \subset L$, $f^{-1}(K)\subset G$ and 
$f^{-1}(L)\subset H$. Since $G$ and $H$ are disjoint, so are $K$ and $L$, and $A \subset \pi g\widehat{D}$-$int(K)$, $B \subset \pi g\widehat{D}$-$int(L)$ and 
$\pi g\widehat{D}$-$int(K) \cap \pi g\widehat{D}$-$int(L)$ = $\phi$. Therefore, $Y$ is softly $\pi g\widehat{D}$-normal.

\begin{corollary}
    If $f: X\longrightarrow Y$ is almost continuous and almost closed surjection and $X$ is a normal space, then $Y$ is softly $\pi g\widehat{D}$-normal.
\end{corollary}
\textbf{Proof.}
    Since every almost closed function is almost $\pi g\widehat{D}$-closed so $Y$ is softly $\pi g\widehat{D}$-normal.

\section{Conclusion}
\label{sec:conclusion}
In this paper, we introduce and study a new weak class of spaces, namely $\pi g\widehat{D}$-normal spaces by using $\pi g\widehat{D}$-open sets. The relationships among normal, strongly quasi-normal, almost softly $\pi g\widehat{D}$- normal spaces and some other topological spaces, also we introduced the concept of $\pi g\widehat{D}$- normal, almost $\pi g\widehat{D}$- normal, quasi $\pi g\widehat{D}$- normal, mildly $\pi g\widehat{D}$- normal. The main result of this paper is that softly $\pi g\widehat{D}$- normality is a topological property and it is a hereditary property with respect to closed domain subspaces. Moreover, we obtain some new characterizations and preservation theorems of softly $\pi g\widehat{D}$- normal spaces. Of course, the entire content will be a successful tool for the researchers for finding the way to obtain the results in the context of such types of normal spaces have many possibilities of applications in digital topology and computer graphics.

\section{Conflict of Interest}
\label{sec:conflict of Intrest}
We certify that this work is original, has never been published before, and isn't being considered for publication anywhere else at this time. This publication is free from any conflicts of interest. As the Corresponding Author, I certify that each of the listed authors has read the paper and given their approval for submission.


\end{document}